\documentclass[12pt,reqno]{amsart}
\usepackage{graphicx,subfigure}
\usepackage{hyperref}
\newcommand{\vertiii}[1]{{\left\vert\kern-0.25ex\left\vert\kern-0.25ex\left\vert #1
    \right\vert\kern-0.25ex\right\vert\kern-0.25ex\right\vert}}
\hypersetup{colorlinks=true, citecolor=blue, linkcolor=red}

\numberwithin{equation}{section} \numberwithin{figure}{section}
\numberwithin{table}{section} 
{ \theoremstyle{plain}
\newtheorem{theorem}{Theorem}[section]

}

\begin{document}

\title[On the nonlinear Schr\"{o}dinger equation in the half-space]
{One-dimensional symmetry of positive bounded solutions to the nonlinear Schr\"{o}dinger equation in the half-space in two and three dimensions}



\author[Christos Sourdis]{Christos Sourdis}% \\ \\ \bigskip Department of Mathematics and Applied Mathematics  \\ University of
%Crete \\  00302106919917} %\protect\footnote{csourdis@tem.uoc.gr}}
\address{Department of Mathematics, National and Kapodistrian University of
Athens, Greece}
\email{sourdis@uoc.gr}

\maketitle

\begin{abstract}
We are concerned with the half-space Dirichlet problem
\[\left\{\begin{array}{ll}
  -\Delta v+v=|v|^{p-1}v & \textrm{in}\ \mathbb{R}^N_+, \\
  v=c\ \textrm{on}\ \partial\mathbb{R}^N_+, &\lim_{x_N\to \infty}v(x',x_N)=0\ \textrm{uniformly in}\ x'\in\mathbb{R}^{N-1},
\end{array}\right.
\]
where $\mathbb{R}^N_+=\{x\in \mathbb{R}^N \ : \ x_N>0\}$ for some $N\geq 2$, and $p>1$, $c>0$ are constants.
It was shown recently  by Fernandez and Weth [Math. Ann. (2021)] that there exists an explicit number $c_p\in (1,\sqrt{e})$, depending only on $p$, such that for $0<c<c_p$ there are infinitely many bounded positive solutions, whereas, for $c>c_p$ there are no bounded positive solutions.
If $N=2, \ 3$, we show that in the case where $c = c_p$
%\[
%v(x',x_N)=c_p\left[\cosh\left(\frac{p-1}{2}x_N\right)
%\right]^{-\frac{2}{p-1}}
%\]
there is no other bounded positive solution besides the one-dimensional one.
\end{abstract}

\section{Introduction}
Recently in \cite{weth}, the authors studied the half-space Dirichlet problem
\begin{equation}\label{eqEq}\left\{\begin{array}{ll}
  -\Delta v+v=|v|^{p-1}v & \textrm{in}\ \mathbb{R}^N_+, \\
 & \\
  v=c\ \textrm{on}\ \partial\mathbb{R}^N_+, &\lim_{x_N\to \infty}v(x',x_N)=0\ \textrm{uniformly in}\ x'\in\mathbb{R}^{N-1},
\end{array}\right.
\end{equation}
where $\mathbb{R}^N_+=\{x\in \mathbb{R}^N \ : \ x_N>0\}$ for some $N\geq 1$, and $p>1$, $c>0$ are constants.
We note that $u(x,t)=e^{it}v(x)$ is a standing wave solution to the focusing nonlinear Schr\"{o}dinger equation with the odd power nonlinearity and exponent $p$.

If $N=1$, then the corresponding ODE has a unique positive even solution that decays to zero at infinity, it is given explicitly by the following formula
\begin{equation}\label{eqw0}
t \rightarrow w_0(t)=c_p\left[\cosh\left(\frac{p-1}{2}t\right)
\right]^{-\frac{2}{p-1}}\ \textrm{with}\ c_p=\left(\frac{p+1}{2} \right)^{\frac{1}{p-1}}=w_0(0)=\sup_{t\in \mathbb{R}}w_0(t).
\end{equation}
Still for $N=1$, it was shown in the aforementioned reference that if $0<c<c_p$ then (\ref{eqEq}) possesses exactly two positive solutions given by
\[
t\to w_0(t+t_{c,p})\ \textrm{and}\ t\to w_0(t-t_{c,p})
\]
with
\[
t_{c,p}=\frac{2}{p-1}\ln\left(\sqrt{\frac{p+1}{2c^{p-1}}}+\sqrt{\frac{p+1}{2c^{p-1}}-1} \right);
\]
 if $c=c_p$ then $w_0$ is the unique positive solution; if $c>c_p$ then there are no positive solutions. We note in passing that the above solutions play an important role in a class of boundary layer problems (see \cite{butuzov}).

If $N\geq 2$, $p>1$, and $0<c<c_p$, using variational methods, it was shown in the same reference that (\ref{eqEq}) admits at least three positive bounded solutions that are geometrically distinct in the sense that they are not translates of each other in the $x'$ direction. In particular, under the further restriction that  $p+1$ is smaller than the critical Sobolev exponent in $\mathbb{R}^N$, $N\geq 2$, then (\ref{eqEq}) admits a positive bounded solution of the form
\[
x\to w_0(x_N+t_{c,p})+u(x) \ \textrm{with}\ u\in H_0^1(\mathbb{R}^N_+)\setminus \{0\} \ \textrm{nonnegative}.
\]

On the other hand, if $c>c_p$, $p>1$, it was shown therein that (\ref{eqEq}) has   no bounded positive solutions. This was accomplished by means of the famous sliding method \cite{berest}.
Roughly speaking, assuming to the contrary that there existed such a solution $v$,  they were able to slide below it  the translates $w_0(x_N+\lambda)$, as $\lambda$ varies continuously starting from very large all the way to $-\infty$. This means that $w_0(x_N+\lambda)\leq v(x)$ in $\mathbb{R}^N_+$ for all $\lambda \in \mathbb{R}$, from where one is led at once to a contradiction. This approach hinges on the translation invariance of the PDE in (\ref{eqEq}) with respect to $x_N$ and the  strong maximum principle applied to the difference of two solutions. The loss of compactness in the $x_N$ direction is compensated by the fact that the linearized PDE at a function that satisfies the uniform decay at infinity enjoys the maximum principle in the half-space $\{x_N>L\}$ with $L\gg1$. On the other side, the loss of compactness in the remaining directions is compensated by the translation invariance property of (\ref{eqEq}) with respect to these (the assumption that $v$ is bounded is used crucially here). The point is that the sliding process cannot get 'stuck' at the boundary of the half-space because there $v=c>c_p=\max w_0$ holds.

 Still in the same reference, it was posed as an interesting open question whether the function $x
\to w_0(x_N)$ is
the unique bounded positive solution to (\ref{eqEq}) in the case $c = c_p$. On the one hand, it seems natural to establish such a uniqueness result also with the help of
a sliding argument, but additional difficulties appear in
the case $c = c_p$ as will be apparent shortly, and non-uniqueness remains a possibility.
On the other hand, one can readily verify that the previously mentioned sliding process starting from $\lambda\gg 1$ can be continued at least  up to $\lambda=0$ (since for  $\lambda>0$ one cannot get stuck on $\partial \mathbb{R}^N_+$, keep in mind that $w_0'(t)<0$ for $t>0$). We point out that
the sliding argument can be initiated exactly as in  Step 1
of the proof of Theorem 6.1 in \cite{weth}; the fact that this process cannot get stuck at some $\lambda_*>0$ (see (6.10) in \cite{weth} for the precise definition) can be proven as in Step 2 of the aforementioned proof (using that on the boundary of the half-space we have $w_0(0+\lambda_*)<w_0(0)=c_p=v(\partial \mathbb{R}^N_+)$). Therefore, with the ordering $w_0\leq v$ at hand, by the strong maximum principle applied in the linear equation for $v-w_0$ (see for instance \cite[Ch. 9]{evans}), we further deduce that
\begin{equation}\label{eqOrder}
  w_0(x_N)<v(x), \ x\in \mathbb{R}^N_+, \ \textrm{or}\ w_0\equiv v.
\end{equation}
In an attempt to exclude the first scenario by pushing further the sliding argument for $\lambda<0$, one has to inevitably  face  the possibility of a sequence $x_j'\in \mathbb{R}^{N-1}$ such that
$\partial_{x_N}v(x_j',0)\to w_0'(0)$ and $|x_j'|\to \infty$ as $j\to \infty$. It is not clear to us how to exclude this behaviour for general $N\geq 2$.

Nevertheless, if $N=2,\ 3$ we are able to exclude the first scenario in (\ref{eqOrder}) by adapting some ideas from the proof of the famous De Giorgi conjecture in the plane (see \cite{bcn,gui}).
In this regard, let us briefly discuss some of the main differences between the two problems. In De Giorgi's conjecture, the solution is assumed to be monotone in some direction. Here, not having a-priori knowledge of such a property, we will instead exploit the ordering (\ref{eqOrder}) and the convexity of the nonlinearity of the problem at hand. In fact, our viewpoint is closer in spirit to that of \cite{dancer,farina}, where the authors establish the  one-dimensional symmetry  of bounded,  stable solutions to semilinear elliptic equations of the form $\Delta v+f(v)=0$ in the whole plane $\mathbb{R}^2$ with $f$ locally Lipschitz (a stronger version of De Giorgi's conjecture).  The fact that we are dealing with the half-space and not the full space also creates some technical difficulties   in applying the approach of the aforementioned references.

More precisely, in this short note we will establish the following one-dimensional symmetry result which answers the aforementioned question of \cite{weth} for $N=2,\ 3$.

\begin{theorem}\label{thm}
  If $N=2,\ 3$ and $p>1$, then the only positive bounded solution of (\ref{eqEq}) with $c=c_p$  is $v(x)=w_0(x_N)$, where $c_p$ and $w_0$ are as in (\ref{eqw0}).
\end{theorem}

The proof of the above result will be given in the following section.

\section{Proof of Theorem \ref{thm}}
\begin{proof}

As we have already explained, the relation (\ref{eqOrder}) is valid. We wish to show that the second alternative is the one which holds. To this end, let us argue by contradiction and suppose
that
\begin{equation}\label{eqordinaryb}
  w_0(x_N)<v(x',x_N), \ (x',x_N)\in \mathbb{R}^N_+.
\end{equation}

Let
\[
z(x',x_N)=w_0'(x_N),\  (x',x_N) \in \overline{\mathbb{R}^N_+}.
\]
Then, $z$ is a negative solution to the linearized problem
\begin{equation}\label{eqEqzzz}
  -\Delta z+z=pw_0^{p-1}(x_N)z \ \textrm{in}\ \mathbb{R}^N_+, \
  z=0\ \textrm{on}\ \partial\mathbb{R}^N_+.
\end{equation}

Let us also consider the difference
\[
\varphi(x',x_N)=v(x',x_N)-w_0(x_N),\ (x',x_N) \in \overline{\mathbb{R}^N_+}.
\]
The main observation is that, owing to the convexity of the nonlinearity $v^p$, $p>1$,  $\varphi$ is a positive supersolution of (\ref{eqEqzzz}) (see also \cite[Prop. 1.3.2]{dupaBook}), i.e.
\begin{equation}\label{eqEqphi}   -\Delta \varphi+\varphi\geq pw_0^{p-1}(x_N)\varphi \ \textrm{in}\ \mathbb{R}^N_+, \
  \varphi=0\ \textrm{on}\ \partial\mathbb{R}^N_+.
\end{equation}
Indeed, since both $v$ and $w_0$ satisfy the same PDE in (\ref{eqEq}), by the mean value theorem we find that
\begin{equation}\label{eqlinear}
  -\Delta \varphi+\varphi=p \theta^{p-1}(x',x_N)\varphi,
\end{equation}
where $w_0(x_N)<\theta(x',x_N)<v(x',x_N)$, $(x',x_N)\in \mathbb{R}^N_+$. The desired relation  (\ref{eqEqphi}) now follows at once from the assumption that $p>1$ and recalling that both $v$ and $w_0$ are equal to $c_p$ on $\{x_N=0\}$.

Armed with the above information,
following \cite{kowas}, we will adapt the method of \cite{bcn,gui} in order to show that
\begin{equation}\label{eqdesire}
  z\equiv a\varphi \ \textrm{for some}\ a\in \mathbb{R}.
\end{equation}
This would imply that $v$ depends only on $x_N$ and therefore has to coincide with $w_0$, which is the unique positive solution of the ODE that is equal to $c_p$ at the origin (this can easily be seen from the phase portrait), in contradiction to (\ref{eqordinaryb}).

The rest of the proof is devoted to showing the validity of (\ref{eqdesire}). To this end, it follows readily from (\ref{eqEqzzz}) and (\ref{eqEqphi}), keeping in mind the signs of $z$, $\varphi$, that
\begin{equation}\label{eqdiv}
  \textrm{div}\left(\varphi^2 \nabla  \left(\frac{z}{\varphi}\right) \right)\leq 0\ \textrm{in}\ \mathbb{R}^N_+.
\end{equation}
For convenience, let us set
\[
g=\frac{z}{\varphi}.
\]
For $R>0$, we let
\[
\xi_R(r)=\zeta\left( \frac{r}{R}\right),\ r=\sqrt{|x'|^2+x_N^2},
\]
where $\zeta$ is some fixed smooth cutoff function such that
\[
\zeta(s)=\left\{\begin{array}{ll}
                  1, & 0\leq s\leq 1, \\
                  0, & s\geq 2,
                \end{array}
 \right.\ \textrm{and}\ \zeta'\leq 0.
\]
We multiply (\ref{eqdiv}) by $g\xi_R^2\leq 0$ and integrate the resulting inequality to get
\[
\int_{x_N>\varepsilon}^{}g\xi_R^2\textrm{div}\left(\varphi^2 \nabla  g \right)dx'dx_N\geq 0\ \textrm{for any}\ \varepsilon>0.
\]
Using the divergence theorem, we obtain
\[
\int_{x_N=\varepsilon}^{}g\xi_R^2\varphi^2 \nabla  g \cdot \nu dx'-\int_{x_N>\varepsilon}^{}\varphi^2\left(\xi_R^2\nabla  g +g\nabla(\xi_R^2)\right)\cdot\nabla g dx'dx_N\geq 0,
\]where $\nu=(0,\cdots,0,-1)$.
Setting
\[\Omega =
\left\{
(x', x_N) \ : \ R<r< 2R\ \textrm{and}\ x_N  > \varepsilon
\right\},\]
it follows that
\begin{equation}\label{eqHuge}
\begin{array}{rcl}
    \int_{x_N>\varepsilon}^{}\varphi^2\xi_R^2|\nabla  g|^2   dx'dx_N & \leq & -\int_{\Omega}^{}\frac{2}{R}g\varphi^2\xi_R\zeta'\left(\frac{r}{R} \right)\nabla r \cdot \nabla g dx'dx_N \\
     &  &  \\
     &  &  +\int_{x_N=\varepsilon}^{}g\xi_R^2\varphi^2 \nabla  g \cdot \nu dx'\\
     &  &  \\
     & \leq &\frac{2}{R}\left\{\int_{\Omega}^{}g^2\varphi^2\left|\zeta'\left(\frac{r}{R} \right)\right|^2dx'dx_N \right\}^\frac{1}{2}\left\{\int_{\Omega}^{}\varphi^2
\xi_R^2|\nabla g|^2dx'dx_N \right\}^\frac{1}{2}  \\
    &  &  \\
     &  & +\int_{x_N=\varepsilon}^{}g\xi_R^2\varphi^2 \nabla  g \cdot \nu dx'.
  \end{array}
\end{equation}
For the last term of the above inequality, we easily see that
\[\begin{array}{rcl}
    \int_{x_N=\varepsilon}^{}g\xi_R^2\varphi^2 \nabla  g \cdot \nu dx' & = & \int_{x_N=\varepsilon}^{}g\xi_R^2 (\varphi\nabla z-z\nabla \varphi) \cdot \nu dx' \\
      &  &  \\
     & = & -\int_{x_N=\varepsilon}^{}g\xi_R^2 (\varphi w_0''(\varepsilon)-z \varphi_{x_N})  dx' \\
     &  &  \\
     & = & -\int_{x_N=\varepsilon}^{}\xi_R^2 \left( w_0'(\varepsilon)w_0''(\varepsilon)-\frac{z^2}{\varphi}   \varphi_{x_N}\right)  dx'.
  \end{array}
\]
By applying Hopf's boundary point lemma in (\ref{eqlinear}) (see for instance \cite[Ch. 9]{evans}), we have
$
 \varphi_{x_N}(x',0)>0$, $x'\in \mathbb{R}^{N-1}$.
Hence,   using the fact that
$z (x',0) = w_0'(0)=0$ and L'Hospital's rule, for any $x' \in \mathbb{R}^{N-1}$, we obtain
\[
\lim_{x_N\to 0^+}\frac{z^2(x',x_N)}{\varphi(x',x_N)}  \varphi_{x_N}(x',x_N)=0.
\]
So, for fixed $R$, we infer that
\[
\int_{x_N=\varepsilon}^{}g\xi_R^2\varphi^2 \nabla  g \cdot \nu dx'\to 0\ \textrm{as}\ \varepsilon\to 0.
\]
Consequently, letting $\varepsilon\to 0$ in (\ref{eqHuge}), recalling the definition of $\zeta$, and since $g\varphi=w_0'$,
we arrive at
\begin{equation}\label{eqRr-}
\begin{split}
\int_{\mathbb{R}^N_+}^{}\varphi^2\xi_R^2|\nabla  g|^2   dx'dx_N\ \ \ \ \ \ \ \ \ \ \ \ \ \ \ \ \ \ \ \ \ \ \ \ \ \ \ \ \ \ \ \ \ \ \ \ \ \ \ \ \ \ \ \ \ \ \ \ \ \ \ \ \ \ \ \ \ \ \ \ \ \ \ \ \ \ \ \ \ \ \ \ &\\ \leq \frac{2}{R}\left\{\int_{\{R<r<2R,\ x_N>0\}}^{} |w_0'(x_N)|^2\left|\zeta'\left(\frac{r}{R} \right)\right|^2dx'dx_N\right\}^\frac{1}{2}\left\{\int_{\{R<r<2R,\ x_N>0\}}^{}\varphi^2\xi_R^2|\nabla  g|^2dx'dx_N\right\}^\frac{1}{2},\ \ \ \ \ \ \ \ &
\end{split}
\end{equation}(note also that, again by L'Hospital's rule, the function $\varphi |\nabla g|$ that appears in the above integrals can actually be extended continuously at $x_N=0$).
Now,  using that $w_0'$ decays exponentially fast as $x_N\to \infty$, we can estimate the first integral of the righthand side  as follows:
\begin{equation}\label{eqRr}
\begin{array}{rcl}
  \int_{\{R<r<2R,\ x_N>0\}}^{} |w_0'(x_N)|^2\left|\zeta'\left(\frac{r}{R} \right)\right|^2dx'dx_N & \leq &    \|\zeta'\|^2_{L^\infty(1,2)}
  \int_{\{|x'|<2R,\ x_N>0\}}^{} |w_0'(x_N)|^2dx' dx_N \\
    &   &   \\
    & = & \|\zeta'\|^2_{L^\infty(1,2)}
  \int_{0}^{\infty} |w_0'(x_N)|^2dx_N \int_{|x'|<2R}^{}dx'\\
&&\\
&= & \|\zeta'\|^2_{L^\infty(1,2)}\|w_0'\|^2_{L^2(0,\infty)}\int_{|x'|<2R}^{}dx'\\
&& \\
&=&  CR^{N-1},
\end{array}
\end{equation}
for some constant $C>0$ that is independent of $R>1$. Hence, by combining (\ref{eqRr-}) and (\ref{eqRr}), we deduce that
\begin{equation}\label{eqRr---}\int_{\mathbb{R}^N_+}^{}\varphi^2\xi_R^2|\nabla  g|^2   dx'dx_N\leq 2\sqrt{C}R^{\frac{N-3}{2}}\left\{\int_{\{R<r<2R,\ x_N>0\}}^{}\varphi^2\xi_R^2|\nabla  g|^2dx'dx_N\right\}^\frac{1}{2}.\end{equation}
   Since $N\leq 3$, the above relation and the monotone convergence theorem imply that
\[
\int_{\mathbb{R}^N_+}^{}\varphi^2|\nabla  g|^2   dx'dx_N<\infty.
\]
In turn, letting $R\to \infty$ in (\ref{eqRr---}), thanks again to the monotone convergence theorem, we conclude that
\[
\int_{\mathbb{R}^N_+}^{}\varphi^2|\nabla  g|^2   dx'dx_N=0,
\]
which gives us the desired identity (\ref{eqdesire}).
\end{proof}

\subsection*{Acknowledgments}The author would like to express his thanks to the anonymous referee for carefully reading the paper and offering some pertinent remarks. Moreover, the author would like to thank IACM of FORTH, where
this paper was written, for the hospitality.
This work has received funding from the Hellenic Foundation for Research and Innovation (HFRI) and the General Secretariat for Research and Technology (GSRT), under grant agreement No 1889.

\end{document}